\documentclass[letterpaper, 10 pt, conference]{IEEEtran}
\IEEEoverridecommandlockouts

\usepackage[numbers,sort&compress]{natbib}

\usepackage[utf8]{inputenc}             \usepackage[T1]{fontenc}                \usepackage[bookmarks=true]{hyperref}   \usepackage{url}                        \usepackage{booktabs}                   \usepackage{colortbl}                   \usepackage{amsfonts}                   \usepackage{nicefrac}                   \usepackage{microtype}                  \usepackage{multicol}

\usepackage{xcolor}         \usepackage{xcolor-material}

\usepackage{graphicx}

\usepackage{dblfloatfix}

\usepackage{mathtools}
\usepackage{framed}
\usepackage{bbm}
\usepackage{xspace}

\usepackage[short]{optidef}       

\usepackage{amssymb}
\usepackage{amsthm}
\usepackage{amsfonts}
\usepackage{amsmath}
\usepackage[capitalise]{cleveref}

\usepackage{bm}
\usepackage{placeins}

\usepackage{booktabs}
\usepackage{multirow}
\usepackage{tabularray}
\UseTblrLibrary{booktabs}

\usepackage{etoolbox}

\usepackage{enumitem}

\makeatletter
\let\MYcaption\@makecaption
\makeatother

\usepackage[font=footnotesize]{subcaption}

\makeatletter
\let\@makecaption\MYcaption
\makeatother

\usepackage{braket}                 \usepackage{xparse}                 

\usepackage[scr=boondox, bb=dsserif]{mathalpha}  \usepackage{bm}

\usepackage{siunitx}
\DeclareSIUnit\ft{ft}               

\usepackage{algorithm}
\usepackage{algpseudocode}
\usepackage{algorithmicx}

\usepackage{nicematrix}             

\usepackage{thmtools}
\declaretheoremstyle[
bodyfont=\normalfont
]{normalstyle}

\usepackage{thm-restate}
\declaretheorem[name=Theorem]{thm}

\declaretheorem[name=Problem]{problem}

\usepackage[framemethod=TikZ]{mdframed}
\mdfdefinestyle{GreyFrame}{nobreak,
    linecolor=white,
    outerlinewidth=0pt,
    roundcorner=2pt,
    leftmargin=-3pt,
    rightmargin=-2pt,
    innertopmargin=4pt, innerbottommargin=2pt, innerrightmargin=3pt,
    innerleftmargin=3pt,
backgroundcolor=MaterialBlueGrey100!20
}

\mdfdefinestyle{ThmFrame}{nobreak,
    linecolor=white,
    outerlinewidth=0pt,
    roundcorner=2pt,
    leftmargin=-3pt,
    rightmargin=-3pt,
    innertopmargin=0pt, innerbottommargin=2pt, innerrightmargin=3pt,
    innerleftmargin=3pt,
    backgroundcolor=MaterialBlueGrey100!15} \definecolor{greentc}{HTML}{52B055}
\definecolor{bluetc}{HTML}{168EF1}
\definecolor{tabblue}{HTML}{1f77b4}
\definecolor{tabred}{HTML}{d62728}

\definecolor{ggRed}{HTML}{E24A33}
\definecolor{ggBlue}{HTML}{348ABD}
\definecolor{ggPurple}{HTML}{988ED5}
\definecolor{ggPurpleStronger}{HTML}{6756d6}
\definecolor{ggGreen}{HTML}{8EBA42}
\definecolor{ggGreenDarker}{HTML}{688C2A}

\definecolor{bndPurple}{HTML}{664d8c}

\definecolor{ciGray}{HTML}{727272}

\newcommand*{\unsafe}{\mathrm{unsf}}
\newcommand*{\safe}{\mathrm{safe}}
\newcommand*{\descent}{\mathrm{desc}}

\makeatletter
\newcommand*{\itemequation}[3][]{\item
  \begingroup
    \refstepcounter{equation}\ifx\\#1\\\else  
      \label{#1}\fi
    \sbox0{#2}\sbox2{$\displaystyle#3\m@th$}\sbox4{\@eqnnum}\dimen@=.5\dimexpr\linewidth-\wd2\relax
\ifcase
        \ifdim\wd0>\dimen@
          \z@
        \else
          \ifdim\wd4>\dimen@
            \z@
          \else 
            \@ne
          \fi 
        \fi
      \@latex@warning{Equation is too large}\fi
    \noindent   
    \rlap{\copy0}\rlap{\hbox to \linewidth{\hfill\copy2\hfill}}\hbox to \linewidth{\hfill\copy4}\hspace{0pt}\endgroup
  \ignorespaces 
}
\makeatother 

\newif\ifarxiv

\newcommand*{\figleft}{{\em (Left)}}
\newcommand*{\figright}{{\em (Right)}}

\makeatletter
\newcommand{\pseudoeqref}[1]{\textup{\tagform@{#1}}}
\makeatother \DeclareFontFamily{U}{BOONDOX-calo}{\skewchar\font=45 }
\DeclareFontShape{U}{BOONDOX-calo}{m}{n}{
  <-> s*[1.05] BOONDOX-r-calo}{}
\DeclareFontShape{U}{BOONDOX-calo}{b}{n}{
  <-> s*[1.05] BOONDOX-b-calo}{}
\DeclareMathAlphabet{\mathcalboondox}{U}{BOONDOX-calo}{m}{n}
\SetMathAlphabet{\mathcalboondox}{bold}{U}{BOONDOX-calo}{b}{n}
\DeclareMathAlphabet{\mathbcalboondox}{U}{BOONDOX-calo}{b}{n}

\DeclarePairedDelimiterX{\norm}[1]{\lVert}{\rVert}{#1}
\DeclarePairedDelimiterX{\abs}[1]{\lvert}{\rvert}{#1}

\newcommand{\Rb}{\mathbb{R}}

\newcommand{\T}{^\mathsf{T}}

\newcommand*{\XSet}{\mathcal{X}}
\newcommand*{\USet}{\mathcal{U}}

\newcommand*{\AvoidSet}{\mathcal{A}}

\newcommand*{\nom}{\mathrm{nom}}

\DeclareDocumentCommand\vectorbold{ s m }{\IfBooleanTF{#1}{\boldsymbol{#2}}{\mathbf{#2}}}

\newcommand*\circled[1]{\tikz[baseline=(char.base)]{
    \node[shape=circle,draw,inner sep=1pt] (char) {#1};}} \def\diffd{\mathrm{d}}

\DeclareDocumentCommand\differential{ o g d() }{ \IfNoValueTF{#2}{
		\IfNoValueTF{#3}
			{\diffd\IfNoValueTF{#1}{}{^{#1}}}
			{\mathinner{\diffd\IfNoValueTF{#1}{}{^{#1}}\argopen(#3\argclose)}}
		}
		{\mathinner{\diffd\IfNoValueTF{#1}{}{^{#1}}#2} \IfNoValueTF{#3}{}{(#3)}}
	}
\DeclareDocumentCommand\dd{}{\differential} 

\DeclareDocumentCommand\derivative{ s o m g d() }
{ \IfBooleanTF{#1}
	{\let\fractype\flatfrac}
	{\let\fractype\frac}
	\IfNoValueTF{#4}
	{
		\IfNoValueTF{#5}
		{\fractype{\diffd \IfNoValueTF{#2}{}{^{#2}}}{\diffd #3\IfNoValueTF{#2}{}{^{#2}}}}
		{\fractype{\diffd \IfNoValueTF{#2}{}{^{#2}}}{\diffd #3\IfNoValueTF{#2}{}{^{#2}}} \argopen(#5\argclose)}
	}
	{\fractype{\diffd \IfNoValueTF{#2}{}{^{#2}} #3}{\diffd #4\IfNoValueTF{#2}{}{^{#2}}}}
}
\DeclareDocumentCommand\dv{}{\derivative}

\usepackage[all=normal,paragraphs=tight,floats=tight,tracking=tight]{savetrees}

\title{\LARGE \bf How to Train Your Neural Control Barrier Function: Learning Safety Filters for Complex Input-Constrained Systems}

\DeclareRobustCommand*{\numauthorrefmark}[1]{\raisebox{0pt}[0pt][0pt]{\textsuperscript{\footnotesize\ensuremath{#1}}}}

\author{\IEEEauthorblockN{Oswin So\numauthorrefmark{1}, Zachary Serlin\numauthorrefmark{2}, Makai Mann\numauthorrefmark{2}, Jake Gonzales\numauthorrefmark{2}, Kwesi Rutledge\numauthorrefmark{1}, Nicholas Roy\numauthorrefmark{1}, Chuchu Fan\numauthorrefmark{1} \thanks{\numauthorrefmark{1} Massachusetts Institute of Technology.}
\thanks{\numauthorrefmark{2} MIT Lincoln Laboratory.}
\thanks{\vspace{-0.75em}\linebreak This work was supported by NASA University Leadership initiative (grants \#80NSSC20M0163 and 80NSSC22M0070), and National Science Foundation CAREER Award (grant \#CCF-2238030). This article solely reflects the opinions and conclusions of its authors and not any NASA entity. \textcopyright 2023 Massachusetts Institute of Technology.}
}}

\begin{document}

\maketitle
\thispagestyle{empty}
\pagestyle{empty}

\begin{abstract}
    Control barrier functions (CBFs) have become popular as a safety filter to guarantee the safety of nonlinear dynamical systems for arbitrary inputs.
However, it is difficult to construct functions that satisfy the CBF constraints for high relative degree systems with input constraints.
To address these challenges, recent work has explored learning CBFs using neural networks via neural CBFs (NCBFs).
However, such methods face difficulties when scaling to higher dimensional systems under input constraints.
    In this work, we first identify challenges that NCBFs face during training.
    Next, to address these challenges, we propose policy neural CBFs (PNCBFs), a method of constructing CBFs by learning the value function of a nominal policy,
    and show that the value function of the maximum-over-time cost is a CBF.
We demonstrate the effectiveness of our method in simulation on a variety of systems ranging from toy linear systems to an F-16 jet with a 16-dimensional state space. Finally, we validate our approach on a two-agent quadcopter system on hardware under tight input constraints.
    The project page can be found at \href{https://mit-realm.github.io/pncbf}{https://mit-realm.github.io/pncbf}.
\end{abstract}

\section{Introduction and Related Works}

Techniques employing control barrier functions (CBFs) are powerful tools for safety-critical control of dynamical systems. 
In particular, CBFs can be used as a safety filter to maintain and certify the safety of any system under arbitrary inputs.
This safety guarantee is crucial in order to give users the needed confidence for greater adoption of robotics in safety-critical domains such as autonomous driving \cite{betz2019autonomous}, surgical robotics \cite{haidegger2019autonomy}, and urban air mobility \cite{NASA_urban}.

Despite their theoretical advantages, constructing CBFs in practice remains difficult.
While it is easy to construct a \textit{candidate} CBF, it is much harder to verify the conditions necessary to enjoy the safety guarantees of a \textit{valid} CBF for systems with input constraints.
Consequently, input constraints are often ignored when using CBFs in practice \cite{xu2018safe,lindemann2018control,wilson2020robotarium,pereira2022decentralized}.

\noindent\textbf{CBFs for High Relative Degree Systems under Input Constraints.}
To address the above challenges with CBFs, recent works try to simplify the construction of valid CBFs for high relative degree systems and input constraints.
In particular, backup CBFs constrain the system to states where a fallback controller can maintain safety \cite{gurriet2020scalable,chen2021backup,squires2018constructive,breeden2021high}. These approaches, however, either require knowledge of an invariant set for the fallback controller which is difficult to compute in itself, or require an appropriate predictive horizon that trades between myopic unsafe behavior and performance.

\noindent\textbf{Neural CBFs. } Recently, learning based approaches have been used to learn neural CBFs (NCBFs) that approximate CBFs using neural networks \cite{dawson2022safeSurvey}, part of a more general trend of learning neural certificates \cite{chang2019neural,deshmukh2019learning,saveriano2019learning,srinivasan2020synthesis,robey2020learning,peruffo2021automated,yang2021iterative,zhao2021learning,lindemann2021learning,cortez2022differentiable}.
Owing to the flexibility of neural networks, NCBFs have been extended to handle parametric uncertainties \cite{dawson2022safe}, obstacles with unknown dynamics \cite{yu2023sequential} and multi-agent control \cite{qin2021learning,zhang2023distributed}. However, many existing NCBF approaches do not consider input constraints. Recent work has examined incorporating input constraints into NCBFs \cite{liu2023safe}, but this approach requires solving a minimax problem that can be brittle to solve in practice.

\noindent\textbf{Reachability Analysis. } Reachability analysis provides a powerful tool for analyzing the safety of dynamical systems. Hamilton-Jacobi (HJ) reachability analysis computes the largest control-invariant set \cite{mitchell2005time} and is often computed using grid-based methods \cite{mitchell2008flexible}. Recent works connect the HJ value function with CBFs \cite{choi2021robust}, providing an alternative to Sum-of-Squares programming for automated synthesis. 
However, the curse of dimensionality limits the practical applicability of grid-based solvers to systems with state-dimension smaller than $5$ \cite{mitchell2008flexible}.

\begin{figure}
    \includegraphics[width=\linewidth]{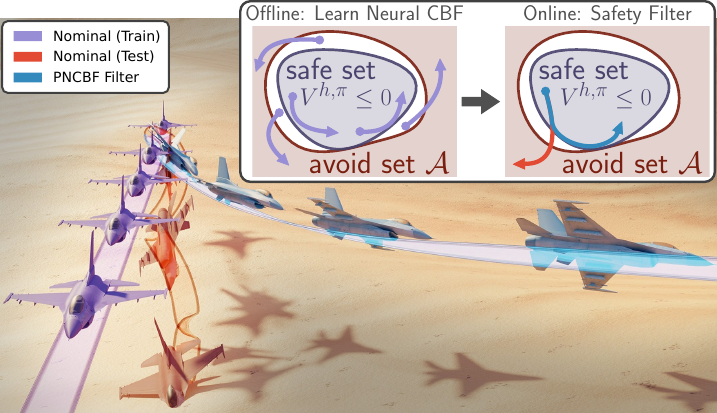}
    \caption{
        We train a Policy Neural CBF (PNCBF) by learning the value function $V^{h,\pi}$ for a given \textcolor{ggPurpleStronger}{\textbf{nominal policy}} offline.
        The sublevel set of $V^{h, \pi}$ contains all states from which the \textcolor{ggPurpleStronger}{\textbf{nominal policy}} remains safe.
        The PNCBF can then be used online as a \textcolor{ggBlue}{\textbf{safety filter}} to ensure the safety of any \textcolor{ggRed}{\textbf{potentially unsafe nominal policy}}.
        Our method avoids the pitfalls of previous Neural CBF approaches and can scale to high dimensional systems such as the F-16.
}
    \vspace{-1.5em}
\end{figure}

\noindent\textbf{Contributions. } We summarize our contributions as follows:
\begin{enumerate}
    \item We identify challenges that existing training methods for Neural CBFs face when under input constraints.
    \item We show that the policy value function is a valid CBF. Using this insight, we propose learning Neural CBFs via the policy value function, thereby bypassing the challenges that previous neural CBF approaches face.
    \item We demonstrate our approach with extensive simulation experiments and show that our method can yield much larger control invariant sets and can scale to higher dimensional systems than current state of the art methods.
    \item We validate our approach on a two-agent quadcopter system on hardware.
\end{enumerate} \section{Preliminaries}
\subsection{Problem Definition}
We consider continuous-time, control-affine dynamics
\begin{equation} \label{eq:dynamics}
    \dot{x} = f(x) + g(x) u,
\end{equation}
where $x \in \XSet \subseteq \Rb^{n_x}, u \in \USet \subseteq \Rb^{n_u}$ and $f, g$ are locally Lipschitz continuous functions. Let $\AvoidSet \subset \XSet$ denote a set of unsafe sets to avoid.
We state the safe controller synthesis problem below.
\begin{mdframed}[style=GreyFrame]
\begin{problem}[Safe Controller Synthesis]\label{prob:safety_problem}
    Given the system \eqref{eq:dynamics} and an avoid set $\AvoidSet \subset \XSet$, find a control policy $\pi : \XSet \to \USet$
    that prevents the system from entering the avoid set $\AvoidSet$, i.e.,
    \begin{equation}
        x_0 \not \in \AvoidSet \implies x_{t} \not \in \AvoidSet, \quad \forall t \geq 0.
    \end{equation}
\end{problem}
\end{mdframed}
Often, we have a nominal policy $\pi_\nom : \XSet \to \USet$ that is performant but may not be safe. In this case, we want our policy $\pi$ to \textit{minimally modify} $\pi_\nom$ to maintain safety.
\begin{mdframed}[style=GreyFrame]
\begin{restatable}[Safety Filter Synthesis]{problem}{}\label{prob:safety_filter_problem}
    Solve \Cref{prob:safety_problem} with the additional desire that $\pi$ is close to $\pi_\nom$. Specifically, we wish to solve the optimization problem
    \begin{mini!}[2]
    {\pi}{ \norm{\pi - \pi_\nom} }{\label{eq:opt:ocp}}{}
    \addConstraint{x_{t}}{\not\in \AvoidSet, \quad\forall t \geq 0,}
    \end{mini!}
    where $\norm{\cdot}$ is some distance metric.
\end{restatable}
\end{mdframed}
In this work, we focus on solving \Cref{prob:safety_filter_problem} with Control Barrier Functions (CBFs), as we describe below.

\subsection{Safety Filter Synthesis with Control Barrier Functions}
We focus on (zeroing) CBFs \cite{wieland2007constructive,xu2015robustness,ames2016control} as a solution to \Cref{prob:safety_filter_problem}. Specifically, let $B : \XSet \to \Rb$ be a continuously differentiable function, $\alpha : \Rb \to \Rb$ be an extended class-$\kappa$ function\footnote{Extended class-$\kappa$ is the set of continuous, strictly increasing functions $\alpha$ such that $\alpha(0) = 0$}, and\footnote{Some works \cite{ames2016control} define the unsafe set to be the zero superlevel set of $B$, while some use the sublevel set. We use the former definition in this work.}
\begin{subequations} \label{eq:cbf:all_constraints}
\begin{align}
    B(x) > 0, \quad \forall x \in \AvoidSet, \label{eq:cbf:classify} \\
    \!\!\!B(x) \leq 0 \implies \inf_{u \in \USet} L_f B(x) + L_g B(x) u \leq -\alpha\big( B(x) \big), \label{eq:cbf:descent}
\end{align}
\end{subequations}
where $L_f B \coloneqq \nabla B\T f$, $L_g B \coloneqq \nabla B\T g$.
Then, $B$ is a CBF, and any control $u$ that satisfies the \textit{descent condition} \eqref{eq:cbf:descent}
renders the sublevel set of $B$ $\Set{x\in\XSet | B(x) \leq 0}$ forward-invariant, i.e., any trajectory starting from within this set remains in this set under such a choice of $u$.
In particular, since \eqref{eq:cbf:descent} is a linear constraint on $u$,
we can solve \Cref{prob:safety_filter_problem} using the following Quadratic Program (QP).
\begin{mini}[2]
{u \in \USet}{ \norm{u - \pi_\nom(x)}^2 \hspace{10em}}
{
\label{eq:cbf:qp}}{}
\addConstraint{L_f B(x) + L_g B(x) u }{\leq -\alpha(B(x))}
\end{mini}

\begin{mdframed}[style=GreyFrame]
\textbf{Challenges with CBF synthesis.}
Define a \textit{candidate} CBF to be any function that satisfies \eqref{eq:cbf:classify}.
If $\mathcal{U}$ is unbounded, since \eqref{eq:cbf:descent} is \textit{linear} in $u$, any candidate CBF $B$ also satisfies \eqref{eq:cbf:descent} if $L_g B \not\equiv 0$. When a system is of high relative degree (i.e., $L_g B(x) \equiv 0$), Higher Order CBFs  (HOCBFs)\cite{nguyen2016exponential,xiao2019control} can be used.

The main challenge to proving that a candidate CBF also satisfies \eqref{eq:cbf:descent} occurs for bounded control sets $\mathcal{U}$ due to actuator limits \cite{breeden2021high}.
Finding a function $B$ such that \eqref{eq:cbf:descent} verifiably holds for \textit{arbitrary} nonlinear dynamics and avoid sets $\mathcal{A}$ is a hard problem that can be solved using Hamilton-Jacobi (HJ) reachability \cite{gurriet2018towards}.
However, HJ reachability is computationally expensive and impractical for systems with more than $5$ dimensions \cite{mitchell2008flexible}.
Consequently, many works that propose CBFs do not consider actuator limits \cite{xu2018safe,lindemann2018control,wilson2020robotarium,pereira2022decentralized}.
\end{mdframed} 

\begin{figure}
    \centering
\includegraphics[width=0.85\linewidth]{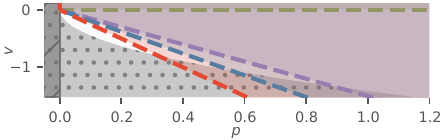}
\caption{\textbf{HOCBF on the double integrator.} On a double integrator with box control constraints $\abs{u} \leq 1$ and the constraint $p \geq 0$, applying different values of the $\alpha$ to the HOCBF candidate $B(x)=-v - \alpha p$ results in different boundaries of the resulting safety filter. However, the only \textit{valid} choice that satisfies the CBF descent condition \eqref{eq:cbf:descent} is $\alpha=0$ (\textcolor{ggGreen}{green}), which disallows any negative velocities and is overly conservative. All other choices of $\alpha$ intersect the true unsafe region (\textcolor{ciGray}{gray dotted}) at some point and violate \eqref{eq:cbf:descent}.}
\label{fig:dbint_hocbf_failure}
\vspace{-0.5em} \end{figure}

One can try to use HOCBFs for automated CBF synthesis of high relative degree systems. As we show next, this can be problematic in the presence of input constraints.
\begin{mdframed}[style=GreyFrame]
\textbf{Challenges of HOCBFs on the Double Integrator.}
Consider the double integrator $\dot{p} = v, \dot{v} = a$, the simplest high relative degree system, with the safety constraint $p \geq 0$.
The HOCBF \textit{candidate} $B(x) = -v - \alpha p$ is a \textit{valid} CBF if and only if $\alpha = 0$ (i.e., disallowing all negative velocities). All other choices of $\alpha$ intersect the true unsafe region and violate \eqref{eq:cbf:descent} (see \cref{fig:dbint_hocbf_failure}).
\end{mdframed}

\subsection{Neural CBFs}

To address the challenges of designing a valid CBF by hand, recent works have proposed learning a CBF using neural networks \cite{dawson2022safe,liu2023safe,yu2023sequential}.

A naive approach to approximating the CBF $B$ with a neural network approximation $B^\theta$ is to encourage satisfying the CBF conditions \eqref{eq:cbf:all_constraints} by minimizing a loss $L$ that penalizes constraint violations over samples of the state space, i.e.,
\begin{subequations}
\begin{align}
    \!\!\!\!\!L_\unsafe{}(\theta, x) &= \left[ \epsilon_{\unsafe} -B_\theta(x) \right]_+, \label{eq:ncbf:loss:unsafe} \\
\!\!\!\!\!L_\descent{}(\theta, x) &= \left[ L_f B_\theta(x) + L_g B_\theta(x) \pi(x) + c B_\theta(x) \right]_+, \label{eq:ncbf:loss:descent} \raisetag{16pt}\\
\!\!\!\!\!L_1(\theta) &= \sum_{\mathclap{x \in \XSet_{\unsafe}}} \; L_\unsafe{}(\theta, x) + \sum_{x \in \XSet} L_\descent{}(\theta, x), \label{eq:ncbf:loss:total}
\end{align}
\end{subequations}
where $\epsilon_{\unsafe}>0$ for the strict inequality in \eqref{eq:cbf:classify},
$\alpha(\cdot)$ is chosen to be linear $x \mapsto c x$ for some $c > 0$, and $\XSet_{\unsafe}$ denotes some superset of $\AvoidSet$. Successful minimization (i.e., zero loss) of \eqref{eq:ncbf:loss:unsafe} implies \eqref{eq:cbf:classify}, and similarly for \eqref{eq:ncbf:loss:descent} and \eqref{eq:cbf:descent}.

However, one problem is that the minimizer of \eqref{eq:ncbf:loss:total} may have a small or even empty forward-invariant set. For example, let $\hat{B}$ be an exponential control-Lyapunov function, i.e.,
\begin{equation}
    \inf_{u \in \USet} L_f \hat{B}(x) + L_g \hat{B}(x) u + \hat{c} \hat{B}(x) \leq 0,  \quad \hat{c} > c.
\end{equation}
Then, $\hat{B} + d$ for all $d>0$ small enough will also have zero loss on \eqref{eq:ncbf:loss:total}. However, the forward-invariant set of $\hat{B} + d$ is the empty set, and hence is not a useful CBF.

\begin{figure}
    \centering
\includegraphics[width=\linewidth]{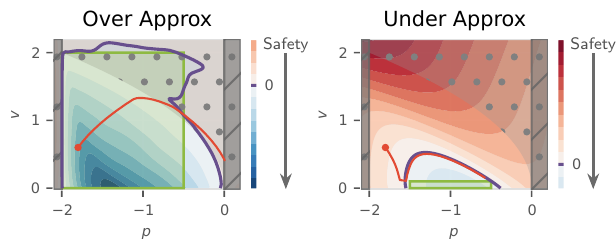}
\caption{\textbf{Over and Underestimation of the Safe Set. } When the safe set (in \textcolor{ggGreenDarker}{green}) is \textit{overestimated} \figleft \,, there are no valid CBFs that can satisfy all the loss terms in \eqref{eq:ncbf:loss:totalsafe} simultaneously. Consequently, the Neural CBF $B^\theta$ has a zero level set (boundary in \textcolor{bndPurple}{purple}) that is larger than the true control-invariant set at the expense of violating the descent condition \eqref{eq:cbf:descent}. In contrast, when the safe set is \textit{underestimated} \figright \,, we can obtain a valid but overly conservative CBF. For comparison, the true unsafe set is shaded in \textcolor{ciGray}{gray} with dots.
}
\label{fig:safeset_overunder} \end{figure}

To address this challenge, many previous works additionally consider a loss term that enforces that $B_\theta \leq 0$ on some safe set $\XSet_{\safe}$ \cite{qin2021learning,dai2022learning,dawson2022safeSurvey,yu2023sequential}, i.e.,
\begin{align}
    L_\safe{}(\theta, x) &= \left[ B_\theta(x) \right]_+, \\
    L_2(\theta) &= \sum_{\mathclap{x \in \XSet_{\safe}}} L_\safe{}(\theta, x) + L_1(\theta). \label{eq:ncbf:loss:totalsafe}
\end{align}
However, the difficulty here is in finding the set $\XSet_{\safe}$. In \cite{qin2021learning,peruffo2021automated}, this is taken to be the set of initial conditions. In \cite{dawson2022safe,dawson2022safeSurvey}, this is assumed to be available, but no details are given for how this set is found in practice. In \cite{yu2023sequential}, this set is
evaluated by rolling out the nominal policy for a fixed number of timesteps. For all these cases, it is not clear whether a valid CBF $B^\theta$ exists such that $B^\theta < 0$ on $\XSet_{\safe}$. The largest-possible $\XSet_{\safe}$ from which a valid CBF can still be found can be obtained using reachability analysis \cite{mitchell2008flexible}. However, the solution of the HJ reachabilty problem yields a CBF directly, rendering the NCBF unnecessary.
Choosing an $\XSet_{\safe}$ that is too large compromises the safety of the resulting CBF, while choosing an $\XSet_{\safe}$ that is too small often results in a forward-invariant set that is too small (see \cref{fig:safeset_overunder}).

An attempt to combat this issue is presented in \cite{liu2023safe}, where a regularization term is added to the loss function to enlarge the sublevel set of the learned CBF $B^\theta$. One issue with this approach is that this regularization term takes a nonzero value for any CBF (including the CBF with the largest zero sublevel set). Hence, the coefficient on this regularization term induces a trade-off between the size of the sublevel set and satisfaction of the CBF constraints \eqref{eq:cbf:classify} and \eqref{eq:cbf:descent}. We provide comparisons against this method in the experiments section \cref{sec:sim_exp}.
 \section{Policy Neural CBFs}
To bypass the above challenges of training a Neural CBF, we now propose policy neural CBFs (PNCBFs), a different approach that does not require knowledge of the safe set but can still recover a large forward-invariant set.
\subsection{Constructing CBFs via Policy Evaluation}
\begin{figure}
    \centering
\includegraphics[width=0.95\linewidth]{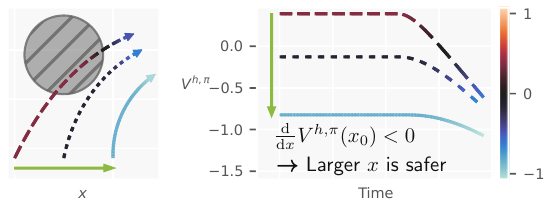}
\caption{\textbf{Understanding the policy value function.}
\figleft \, Trajectories from a nominal policy $\pi$ started from three values of $x_0$. \figright \, The corresponding policy value functions $V^{h,\pi}$ along each trajectory.
$V^{h,\pi}$ is non-increasing along (any) trajectory of $\pi$ and is a CBF. Hence, the gradients of $V^{h,\pi}$ inform the CBF-QP on how to improve safety using the ``knowledge'' of the $\pi$.
}
\label{fig:pncbf_intuition}
\vspace{-1.1em} \end{figure}

We assume that the avoid set $\mathcal{A}$ can be described as the superlevel set of some continuous function $h : \XSet \to \Rb$, i.e.,
\begin{equation} \label{eq:A_def}
    \mathcal{A} = \Set{x \in \XSet | h(x) > 0}.
\end{equation}
Let $\pi : \XSet \to \USet$ be an arbitrary policy, and let $x^\pi_t$ denote the resulting state at time $t$ following $\pi$.
Consider the following \textit{maximum-over-time} cost function
\begin{equation} \label{eq:max_over_time_cost}
    V^{h,\pi}(x_0) \coloneqq \sup_{t \geq 0} h(x^\pi_t).
\end{equation}
It can be shown that $V^{h,\pi}$ satisfies the following Hamilton-Jacobi PDE in the viscosity sense \cite{altarovici2013general}.
\begin{equation} \label{eq:V:pde}
    \max\Big\{ h(x) - V^{h,\pi}(x),\; \nabla V^{h, \pi}(x)\T \left( f(x) + g(x) \pi(x) \right) \Big\} = 0.
\end{equation}
This immediately gives us the following two inequalities
\begin{subequations}\label{eq:V_ineqs}
\begin{align}
    V^{h,\pi}(x) &\geq h(x), \label{eq:V:ineq_h}\\
    \nabla V^{h, \pi}(x)\T \left( f(x) + g(x) \pi(x) \right) &\leq 0, \label{eq:V:ineq_Vdot}
\end{align}
\end{subequations}
from which we have the following theorem.
\begin{mdframed}[style=ThmFrame]
\begin{thm}[Policy value function is a CBF]
    The policy value function $V^{h,\pi}$ is a CBF for \eqref{eq:dynamics} for any $\pi$ and $\alpha>0$.
\end{thm}
\end{mdframed}
\begin{proof}
\eqref{eq:V:ineq_h} and \eqref{eq:A_def} implies \eqref{eq:cbf:classify}. Next, \eqref{eq:V:ineq_Vdot} implies \eqref{eq:cbf:descent} for \textit{any} choice of $\alpha$, since $V(x) \leq 0$ implies that
\begin{equation*}
\pushQED{\qed} 
    \nabla V^{h, \pi}(x)\T \left( f(x) + g(x) \pi(x) \right) \leq 0 \leq -\alpha( V(x) ). \qedhere
\end{equation*}
\end{proof}
Intuitively, the policy value function $V^{h, \pi}$ gives us an upper-bound on the worst constraint violation $h$ in the future under the optimal policy, since using $\pi$ guarantees that $h$ will be at most $V^{h, \pi}$, and the optimal policy will do no worse. Moreover, by following the negative gradient of $V^{h,\pi}$, we can move to states where following $\pi$ leads to a lower maximum value of $h$, i.e., safer states (see \cref{fig:pncbf_intuition}).

Consequently, this provides us with a method to construct CBFs via policy evaluation of \textit{any} policy $\pi$. To make this more concrete, consider the dynamic-programming form of \eqref{eq:max_over_time_cost}:
\begin{equation} \label{eq:avoid_dp}
    V^{h,\pi}(x_0) = \max\left\{ \sup_{0 \leq s \leq t} h(x_s),\enspace V^{h, \pi}(x_t)\right\}.
\end{equation}
Given a nominal policy $\pi$, we can collect rollouts of the system and store tuples $\left( x_0, \max_{0 \leq t \leq T} h(x_t), x_T \right)$. We then minimize the policy evaluation loss on a neural network approximation of the policy value function $V^{h, \pi, \theta}$
\begin{equation} \label{eq:pncbf:loss}
    \!\!\!\!\!L = \norm*{ V^{h,\pi}_\theta(x_t) - \max\left\{ \max_{t \leq s \leq T} h(x_T), V^{h, \pi}_\theta h(x_T) \right\} }^2.
\end{equation}
We summarize the above for training PNCBFs in \Cref{alg:pncbf}.

\noindent After training, we can use $V^{h,\pi}_\theta$ via the CBF-QP \eqref{eq:cbf:qp} to minimally modify the (unsafe) nominal policy to maintain safety.

\begin{algorithm}[t]
    \small
    \caption{Policy Neural CBF}
    \label{alg:pncbf}
    \begin{algorithmic}[1]
        \State{\textbf{input: }Nominal Policy $\pi$}
        \State Collect dataset of tuples $\big( x_0,\, \max_{0 \leq t \leq T} h(x_t),\, x_T \big)$
        \While{not converged}
            \State Minimize loss \eqref{eq:pncbf:loss} over samples from the dataset
        \EndWhile
    \end{algorithmic}
\end{algorithm} 
\vspace{0.25em}

\noindent\textbf{Viewing policy CBFs as policy distillation. } One can interpret policy value functions as policy distillation. More specifically, when $V^{h, \pi}$ is used as a safety filter in the CBF-QP \eqref{eq:cbf:qp} with any \textit{new} nominal policy $\tilde{\pi}$, the forward-invariant set of the resulting CBF-QP controller will be no smaller than that of the original nominal policy $\pi$, as we show next.
\begin{mdframed}[style=ThmFrame]
\begin{thm} \label{thm:pncbf_fwdinvariance}
    Let $V^{h,\pi}$ be a policy value function and let $\tilde{\pi}$ be some other policy. Then, the forward-invariant set under CBF-QP with $V^{h,\pi}$ and $\tilde{\pi}$ is a superset of the forward-invariant set under $\pi$.
\end{thm}
\end{mdframed}
\begin{proof}
    The forward-invariant set under $\pi$ is exactly the zero sublevel set of $V^{h,\pi}$ $\Set{x | V^{h, \pi} \leq 0}$. Since $V^{h,\pi}$ is a CBF, the CBF-QP controller will render this set forward-invariant under any \textit{new} nominal policy $\tilde{\pi}$.
\end{proof}

\begin{mdframed}[style=GreyFrame]
\textbf{Relationship with Hamilton-Jacobi Reachability.}
The policy CBF is also closely related to HJ reachability. As noted in \cite{choi2021robust,tonkens2022refining,tonkens2023patching}, the (optimal) HJ value function is a CBF. This is equivalent to the policy CBF \eqref{eq:max_over_time_cost} with the optimal policy $\pi^*$. The policy CBF can thus be seen as a \textit{relaxation} of optimality that remains a CBF.
\end{mdframed}
For neural networks, policy evaluation can be more attractive than optimization, which requires techniques such as deep reinforcement learning (e.g., \cite{fisac2019bridging,hsu2021safety}) that can be more unstable and computationally expensive.
However, as a middle ground, we next show how policy iteration can be applied to PNCBFs to achieve fast convergence without resorting to a full deep reinforcement learning setup.

\begin{figure*}[t]
    \includegraphics[width=\linewidth]{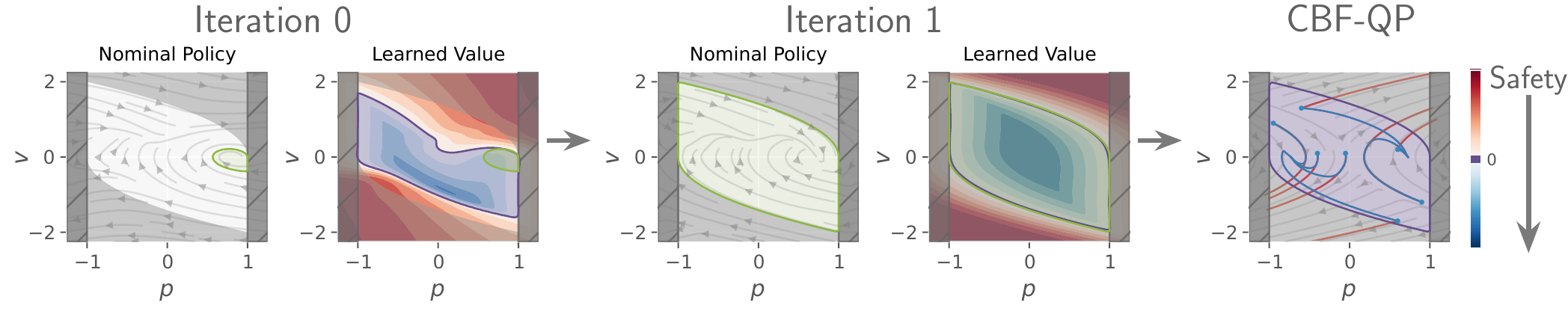}
    \caption{\textbf{Policy iteration on a double integrator.} Starting with a suboptimal nominal policy $\pi$, we learn the value function $V^{h, \pi}$. By treating the CBF-QP of the learned $V^{h, \pi}$ as a new nominal policy, we can repeat this process to perform \textit{policy iteration}. Here, only two iterations are needed to obtain a CBF $V^{h, \pi}$ that almost covers the true control-invariant set. The final CBF-QP controller maintains safety (\textcolor{ggBlue}{blue line}) under any potentially unsafe nominal policy (\textcolor{ggRed}{red line}).}
    \label{fig:policy_iteration}
    \vspace{-0.9em}
\end{figure*}

\subsection{Policy Iteration with PNCBFs}

The choice of the nominal policy $\pi$ is crucial.
In light of \Cref{thm:pncbf_fwdinvariance}, the forward invariant set of the resulting PNCBF controller (via the CBF-QP \eqref{eq:cbf:qp}) is only guaranteed to be no smaller than that of $\pi$. Hence, a poor choice of $\pi$ can result in a small forward-invariant set, resulting in a poor CBF.

To resolve this problem, we use the insight that the policy value function \eqref{eq:max_over_time_cost} is also a (shifted) Lyapunov function. Hence, when using $V^{h, \pi}$ with the CBF-QP \eqref{eq:cbf:qp}, we can hope that the resulting forward-invariant set will be larger than the original policy $\pi$. Nevertheless, this new controller will be no worse than the original policy $\pi$. Thus, we propose to take the PNCBF controller as the \textit{new} nominal policy $\pi^+$ to train a new PNCBF, and iterate this procedure (see \cref{fig:policy_iteration}). 
By treating the application of a CBF-QP as an analytical \textit{policy improvement} and the computation of $V^{h,\pi}$ as policy evaluation, we can interpret this procedure as \textit{policy iteration}, which has been studied extensively for the normal sum-over-time cost structure \cite{bellman1966dynamic} where it enjoys guaranteed convergence at a superlinear rate under certain assumptions \cite{santos2004convergence}.
While it is not clear if this convergence result holds for the maximum-over-time cost structure, we empirically observe fast convergence in only a few iterations, as we show in \cref{sec:sim:pol_iter}.
Also, we observe that using a policy value function with a non-zero discount factor can help with convergence when $\pi$ is far from optimal. We leave an analysis of the interaction between the discount factor and convergence rates to future work. 
\subsection{Discounting and Contraction}
One problem with using \eqref{eq:pncbf:loss} directly as a loss function is that there are undesireable solutions that satisfy this recursive equation. For example, $V^{h,\pi}(x) = a$ for $a$ large enough minimizes \eqref{eq:pncbf:loss}, but is clearly not a solution to \eqref{eq:max_over_time_cost}. 
This is similar to the case in Markov Decision Processes where the \textit{undiscounted} value iteration is not contractive \cite{federgruen1978contraction}.
Hence, instead of \eqref{eq:max_over_time_cost},
we consider the following \textit{discounted} cost, for $\lambda \geq 0$.
\begin{align} \label{eq:max_over_time_disc_cost}
    V^{h,\pi}_\lambda(x_0) &\coloneqq \sup_{t \geq 0} \left\{ \tilde{h}(x_t, \lambda) + e^{-\lambda t} h(x_t) \right\}, \\
\tilde{h}(x_t, \lambda) &\coloneqq \int_0^t \lambda e^{-\lambda s} h(x_s) \diffd{s}.
\end{align}
Taking $\lambda = 0$ recovers the undiscounted problem \eqref{eq:max_over_time_cost}, while $\lambda \to \infty$ yields the solution $V^{h, \pi}_\infty = h$. Hence, different choices of $\lambda$ can be seen as \textit{implicitly} choosing the horizon considered for safety. Similar to \eqref{eq:V:pde}, it can also be shown that $V^{h,\pi}_\lambda$ satisfies the following Hamilton-Jacobi PDE in the viscosity sense (suppressing arguments for conciseness) \cite{altarovici2013general}:
\begin{equation} \label{eq:V:disc_pde}
    \max\Big\{ h - V^{h,\pi}_\lambda,\; \nabla {V^{h, \pi}_\lambda}\T \left( f + g \pi \right) - \lambda (V^{h, \pi}_\lambda - h)\Big\} = 0,
\end{equation}
as well as the following dynamic programming equation
\begin{equation} \label{eq:pncbf:disc_loss}
    V^{h,\pi}_\lambda(x_0) = \max\left\{ \sup_{0 \leq s \leq t} \tilde{h}(x_s, \lambda), \tilde{h}(x_t, \lambda) + e^{-\lambda t} V^{h, \pi}_\lambda(x_t)\right\}.
\end{equation}
While solutions to the PDE \eqref{eq:V:disc_pde} no longer satisfy the CBF constraint \eqref{eq:cbf:descent} for $\lambda > 0$, they do prevent the constant solution from being a minimizer of the corresponding discounted loss. Hence, in practice, we use \eqref{eq:pncbf:disc_loss} instead of \eqref{eq:pncbf:loss}. We start with a small value of $\lambda$ to avoid premature convergence to the constant solution, and gradually decrease it to $0$ as training progresses. 
\begin{mdframed}[style=GreyFrame]
\noindent\textbf{Verification of PNCBF. }
We stress that, without verifying that the learned PNCBF satisfies the descent condition \eqref{eq:cbf:descent}, we can not claim that the PNCBF satisfies the CBF conditions nor claim any safety guarantees.
Verification of NCBFs can be performed using neural-network verification tools \cite{liu2021algorithms}, sampling \cite{bobiti2018automated} or a generalization error bound \cite{qin2021learning}.
\end{mdframed}

Nevertheless, as we show next, empirical results show that our proposed method vastly improves the volume of both the forward-invariant set and the set where the safety filter is permissible to nominal controls compared to baseline methods, including an (unverified) HOCBF candidate.
 \section{Simulation Experiments} \label{sec:sim_exp}

To study the performance of PNCBFs, we perform a series of simulation experiments on high relative degree systems under box control constraints.
We first investigate the qualitative behavior of PNCBFs on simple low-dimensional systems to gain insight into the behavior of the method.
Next, we demonstrate the scalability of PNCBFs to high-dimensional systems by applying it to a ground collision avoidance problem on the F-16 fighter jet \cite{heidlauf2018verification,stevens1992aircraft}.
Finally, we study the behavior of policy iteration on a two-agent quadcopter system to demonstrate the ability of PNCBFs to improve the safety of an initially unsafe nominal policy.

\noindent\textbf{Baselines. } We compare against the following safety filters.
\begin{itemize}
\item \textbf{Neural CBF (NCBF) \cite{dawson2022safe,dawson2022safeSurvey}: } Learning a Neural CBF using \eqref{eq:ncbf:loss:totalsafe}. We choose the safe set to be the set containing the equilibrium point under the nominal policy $\pi$.
\item \textbf{Non-Saturating Neural CBF (NSCBF) \cite{liu2023safe}: } A recent approach that explicitly tackles the problem of input constraints for CBFs by learning a Neural CBF. However, instead of enforcing the derivative condition \eqref{eq:cbf:descent} over the entire state-space as in \cite{dawson2022safe}, this is only enforced on the boundary as in barrier certificates \cite{prajna2004safety}.
\item \textbf{Handcrafted Candidate CBF (CBF) \cite{nguyen2016exponential,xiao2019control}: } We construct a \textit{candidate} CBF via a Higher-Order CBF on $h$ without considering input constraints.
\item \textbf{Approximate MPC-based Predictive Safety Filter (MPC) \cite{wabersich2021predictive}: }
A trajectory optimization problem is solved, imposing the safety constraints while penalizing deviations from the nominal policy.
    We do not assume access to a known forward-invariant set and hence do not impose this terminal constraint.
\item \textbf{Sum-of-Squares Synthesis (SOS) \cite{zhao2023convex}: } When the dynamics are polynomial, a sequence of convex optimization problem can be solved to construct a CBF.
\end{itemize}

All neural networks are trained until convergence and use the same architecture (3 layers of $256$ neurons with $\tanh$ activations).
For PNCBFs, we perform at most $3$ iterations of policy iteration.

\begin{figure*}
    \centering
\includegraphics[width=\linewidth]{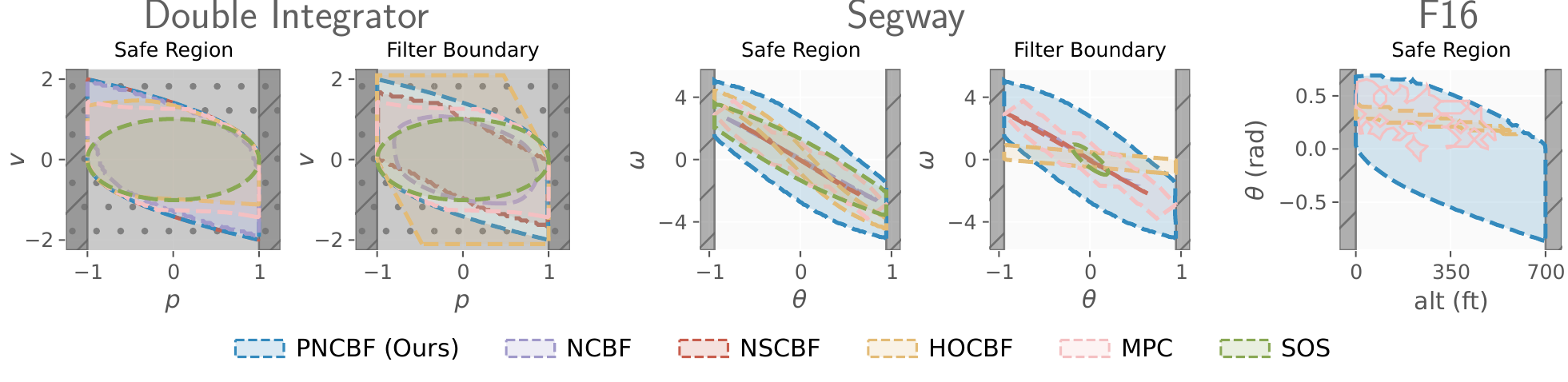}

\vspace{-0.5em}
\caption{\textbf{Safe set and filter boundary on the double integrator, Segway, and F16} 
We plot the initial states from where the safety filter can preserve safety (Safe Region), and states where the nominal policy can influence the output of the safety filter (Filter Boundary). The true unsafe region is shaded in \textcolor{ciGray}{gray dots}.
On the double integrator, ours is the only method that is both maximally safe and permissive. For more complex systems, the performance gap between our method and baseline methods becomes more pronounced, showcasing the benefit of PNCBFs on high-dimensional nonlinear systems.
}
\label{fig:sim_ci}
\vspace{-0.5em} \end{figure*}
\begin{figure*}
    \centering
\includegraphics[width=0.2\linewidth]{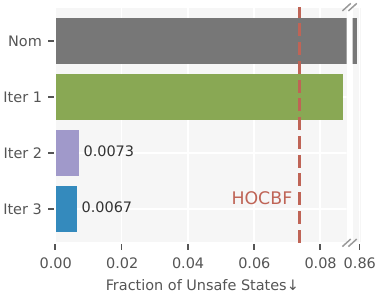}\hfill
\includegraphics[width=0.77\linewidth]{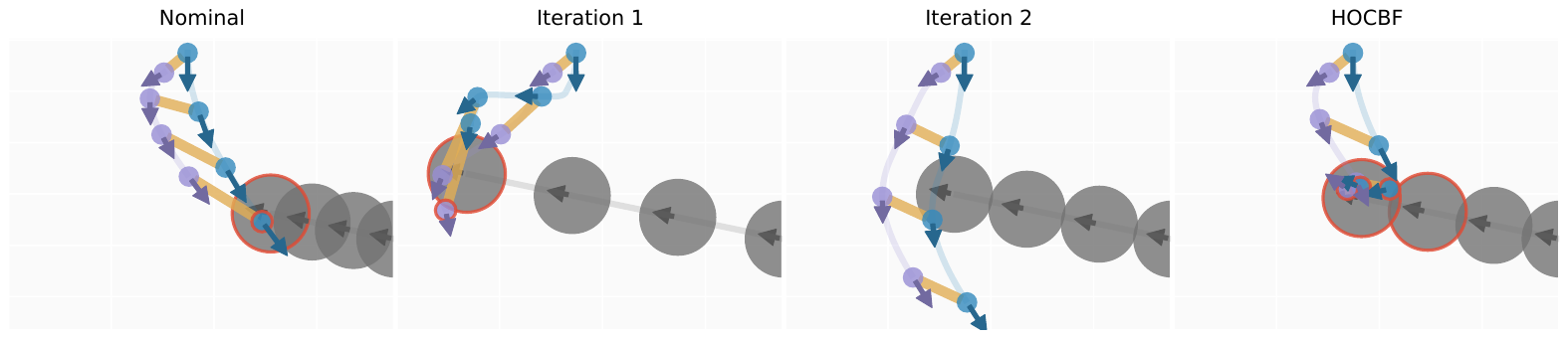}
\caption{\textbf{Policy Iteration on a Two-Agent Quadcopter System.}
\figleft \, In only three iterations, we achieve the smallest volume of unsafe states compared to baseline methods and greatly improving the safety of the original unsafe nominal policy.
Since we may sample states in the true unsafe region, the optimal safety filter will not be able to achieve safety for all sampled states. \figright \, An example of an initial state from which only our method is able to filter out unsafe controls from the nominal policy and prevent a collision (highlighted in \textcolor{ggRed}{red}) with the moving obstacle. 
}
\label{fig:quad_stats}
\vspace{-1.2em} \end{figure*}

\subsection{Qualitative Behavior on a Double Integrator, Segway}
We first perform a general comparison between the different methods on a double integrator and a Segway, two simple systems that can be easily visualized. On the double integrator, safety is defined via position bounds ($\abs{p} \leq 1$), while the Segway asks for the handlebars to stay upright ($\abs{\theta} \leq 0.3\pi$) while remaining within position bounds ($\abs{p} \leq 2.0$). During testing, we use a different nominal policy (i.e., zero control) than the one used during training for PNCBFs.

We visualize the results in \cref{fig:sim_ci},
plotting the region of the state space from where the safety filter preserves safety (Safe Region) and where the nominal policy can influence the output of the safety filter (Filter Boundary). For CBF-based methods, the filter boundary corresponds to the zero level set of the CBF.
On the double integrator, all methods induce forward-invariance on some region of the state space but only our method is both maximally safe and permissive. This trend is even more pronounced on the Segway, where our method is able to find a significantly larger safe set and filter boundary.

\subsection{Scalability to high-dimensional systems with F-16}
Next, we explore the scalability of PNCBF to high dimensional systems. We consider a ground collision avoidance example involving the F-16 fighter jet \cite{heidlauf2018verification,stevens1992aircraft}. Since this system is not control-affine in the throttle, we leave the throttle as the output of a P controller, resulting in a $16$-dimensional state space and a $3$-dimensional control space.
We define safety as a box constraint on the aircraft's altitude.
During testing, we apply an adversarial nominal policy that commands the aircraft to dive nose-first into the ground.

We visualize the results in \cref{fig:sim_ci},
showing a $2$D slice of the state space. Even on a $16$-dimensional state space, we observe that PNCBFs are able to recover a significantly larger region of the safe set compared to other baseline methods.

\subsection{Performance of Policy Iteration} \label{sec:sim:pol_iter}
Finally, to investigate the ability of PNCBFs to learn a safe and permissible safety filter from an initially unsafe nominal policy, we consider a two-agent quadcopter system with a $12$-dimensional state and $4$-dimensional control space that must stay within communication radius while avoiding collisions with a dynamic obstacle. We model each quadcopter as a double integrator with a velocity tracking controller.
The obstacle is assumed to move with constant velocity and direction. The nominal policy moves each quadcopter anticlockwise around a circle, ignoring all constraints.
The obstacle can achieve higher velocities than each quadcopter. Additionally, the velocity tracking controller has a slow response time. Hence, the quadcopters must react well in advance to avoid collisions with the obstacle while staying within communication radius, resulting in a problem with complex safety constraints despite the simple dynamics.

We visualize the results in \cref{fig:quad_stats}.
Although the nominal policy has a high unsafe fraction, policy iteration is able to significantly reduce the unsafe fraction to near $0$ in only two iterations, representing a $90\%$ reduction in unsafe states compared to the next best method.
 \section{Hardware Experiments}
\begin{figure}[t]
    \centering
\includegraphics[width=0.85\linewidth]{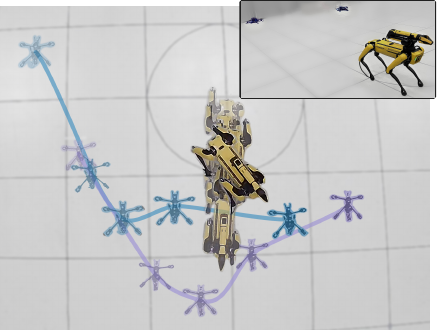}
\caption{\textbf{Two-agent Quadrotor System with Moving Quadruped Obstacle.} Snapshots from a hardware realization of the setup from \cref{sec:sim:pol_iter}.
}
\label{fig:hw_img}
\vspace{-1em} \end{figure}

We further validate our approach in a two-agent quadrotor hardware experiment mirroring the setup from \cref{sec:sim:pol_iter}.
We use two custom drones and use Boston Dynamics’s Spot as a dynamic obstacle.
Velocity setpoints are sent to the drones through the PX4 flight stack.
We visualize the results in \cref{fig:hw_img}, where the PNCBF filters the drone's unsafe nominal policy to avoid collisions with Spot while remaining within communication radius.
For more details, see the supplemental video.

 \section{Conclusion}
By learning the policy value function, we are able to learn Neural CBFs for high relative degree systems under input constraints. Extensive simulation experiments show that our method can yield much larger forward invariant sets and can be deployed on hardware.
One limitation of our method is that it requires an accurate dynamics model. Model errors may cause the learned safety filter to be unsafe on the real system. While this was not a major issue in our hardware experiments, we plan to investigate methods to improve robustness to model errors in future work, such as by learning \textit{robust} CBFs as in \cite{dawson2022safe}.

\bibliographystyle{IEEEtran}
\bibliography{IEEEabrv,references.bib}

\newpage
\onecolumn

\begin{appendices}
\numberwithin{equation}{section}
\renewcommand{\theequation}{\thesection.\arabic{equation}}
\renewcommand{\thesubsection}{\thesection\arabic{subsection}}
\renewcommand{\thesubsectiondis}{\thesubsection}

\crefalias{section}{appendix}
\crefalias{subsection}{appendix}

\setlength{\parindent}{0pt}
\setlength{\parskip}{2pt plus1pt minus0.5pt}
\section{Derivation of the Discounted Hamilton Jacobi PDE} \label{app:sec:hj_derivation}
Define the following \textit{discounted} objective function $V^{h,\pi}_\lambda$ for $\lambda \geq 0$, as in \eqref{eq:max_over_time_disc_cost}:
\begin{align}
    V^{h,\pi}_\lambda(x_0) &\coloneqq \sup_{t \geq 0} \left\{ \tilde{h}(x_t, \lambda) + e^{-\lambda t} h(x_t) \right\}, \label{eq:app:max_over_time_disc_cost} \\
\tilde{h}(x_t, \lambda) &\coloneqq \int_0^t \lambda e^{-\lambda s} h(x_s) \diffd{s}.
\end{align}
We now derive a recursive relationship for $V^{h, \pi}_\lambda$. Starting from \eqref{eq:app:max_over_time_disc_cost} and splitting the supremum into two parts, we have
\begin{equation} \label{eq:app:V:split_sup}
    V^{h,\pi}_\lambda(x_0) 
    = \max\bigg\{ \sup_{0 \leq r \leq t} \left\{ \tilde{h}(x_r, \lambda) + e^{-\lambda r} h(x_r) \right\},\; \underbrace{\sup_{r \geq t} \left\{ \tilde{h}(x_r, \lambda) + e^{-\lambda r} h(x_r) \right\}}_{\circled{\footnotesize$\star$}} \bigg\}
\end{equation}
We next simplify $\circled{\footnotesize$\star$}$ as
\begin{align}
    \circled{\footnotesize$\star$}
    &\coloneqq \sup_{r \geq t} \left\{ \tilde{h}(x_r, \lambda) + e^{-\lambda r} h(x_r) \right\} \\
    &= \int_0^t \lambda e^{-\lambda s} h(x_s) \dd{s} + \sup_{\Delta \geq 0} \left\{ \int_t^{t+\Delta} \lambda e^{-\lambda s} h(x_s) \dd{s} + e^{-\lambda (t+\Delta)} h(x_{t + \Delta}) \right\} \\
    &= \tilde{h}(x_t,\lambda) + e^{-\lambda t} \sup_{\Delta \geq 0} \left\{ \int_0^{\Delta} \lambda e^{-\lambda s} h(x_{t+s}) \dd{s} + e^{-\lambda \Delta} h(x_{t + \Delta}) \right\} \\
    &= \tilde{h}(x_t,\lambda) + e^{-\lambda t} V^{h, \pi}_\lambda(x_t).
\end{align}
Plugging this back into \eqref{eq:app:V:split_sup}, we thus obtain
\begin{equation} \label{eq:app:recursive}
    V^{h,\pi}_\lambda(x_0) 
    = \max\bigg\{ \underbrace{\sup_{0 \leq r \leq t} \left\{ \tilde{h}(x_r, \lambda) + e^{-\lambda r} h(x_r) \right\}}_{\circled{\footnotesize 1}},\; \underbrace{\tilde{h}(x_t, \lambda) + e^{-\lambda t} V^{h, \pi}_\lambda(x_t)}_{\circled{\footnotesize 2}} \bigg\}.
\end{equation}
Finally, we now state and prove the following lemma.
\begin{mdframed}[style=ThmFrame]
\begin{thm} \label{app:thm:hj_pde_disc}
    The function $V^{h,\pi}_\lambda$ is the viscosity solution of the following HJ inequality.
    \begin{equation} \label{eq:app:hj_ineq}
        \max\Big\{ h(x) - V^{h,\pi}_\lambda(x),\; \nabla {V^{h, \pi}_\lambda}(x)\T \big( f(x) + g(x) \pi(x) \big) - \lambda (V^{h, \pi}_\lambda(x) - h(x))\Big\} = 0.
    \end{equation}
\end{thm}
\end{mdframed}
\begin{proof}
    To prove that $V^{h,\pi}_\lambda$ is a viscosity solution, we need to prove that it is both a super-solution and a sub-solution. We first prove the super-solution property.

\textbf{Viscosity super-solution.}
    From \eqref{eq:app:recursive}, since $V^{h,\pi}_\lambda \geq \circled{\footnotesize1}$ by definition of $\max$,
    \begin{equation}
        V^{h,\pi}_\lambda(x_0) \geq \sup_{0 \leq r \leq t} \left\{ \tilde{h}(x_r, \lambda) + e^{-\lambda r} h(x_r) \right\} \geq h(x_0).
    \end{equation}
    Hence, 
    \begin{equation}
      h(x_0) - V^{h,\pi}_\lambda(x_0) \leq 0 \label{eq:app:supersol:p1}.
    \end{equation}
    Next, we use $V^{h,\pi}_\lambda \geq \circled{\footnotesize2}$ to obtain that, for any $t \geq 0$,
    \begin{equation}
        \int_0^0 \lambda e^{-\lambda s} h(x_s) \dd{s} + e^{-\lambda 0} V^{h,\pi}_\lambda(x_0) = V^{h,\pi}_\lambda(x_0) \geq \int_0^t \lambda e^{-\lambda s} h(x_s) \dd{s} + e^{-\lambda t} V^{h,\pi}_\lambda(x_t).
    \end{equation}
    By applying classical arguments in viscosity theory, we get
    \begin{equation}
        \dv{t} \left[
        \lambda \int_0^t e^{-\lambda \tau} h(x_\tau) \,\diffd{\tau} + e^{-\lambda t} V(x_t)
        \right]_{t=0}
        = V_x\T \big( f(x_t) + g(x_t) \pi(x_t) \big) - \lambda( V(x_t) - h(x_t))
        \leq 0 \label{eq:app:supersol:p2}.
    \end{equation}
    Combining \eqref{eq:app:supersol:p1} and \eqref{eq:app:supersol:p2}, we obtain that
    \begin{equation}
        \max\Big\{ h(x) - V^{h,\pi}_\lambda(x),\; \nabla {V^{h, \pi}_\lambda}(x)\T \big( f(x) + g(x) \pi(x) \big) - \lambda (V^{h, \pi}_\lambda(x) - h(x))\Big\} \leq 0.
    \end{equation}
    In other words, $V^{h,\pi}_\lambda$ is a viscosity supersolution.

\textbf{Viscosity sub-solution. }
    We prove this by considering two cases.

    \underline{Case 1}: Suppose that
    \begin{equation}
        h(x) \geq V^{h,\pi}_\lambda(x).
    \end{equation}
    Then,
    \begin{equation} \label{eq:app:subsol:p1}
        \max\Big\{ h(x) - V^{h,\pi}_{\lambda}(x), \; \nabla {V^{h,\pi}_{\lambda}}\T \big( f(x) + g(x) \pi(x) \big) - \lambda( V^{h,\pi}_{\lambda}(x) - h(x) ) \Big\} \geq 0.
    \end{equation}

    \underline{Case 2}: On the other hand, suppose that
    \begin{equation}
        h(x_0) < V^{h,\pi}_{\lambda}(x_0).
    \end{equation}
    By the continuity of $h$ and $V^{h,\pi}_{\lambda}$, there exists some upper bound $t > 0$ such that for all $r \in [0, t]$, we have that $h(x_r) < V^{h,\pi}_\lambda(x_r)$. Multiplying by $e^{-\lambda r}$ and adding $\tilde{h}(x_r, \lambda)$ to both sides, we get
    \begin{equation}
        \int_0^r \lambda e^{-\lambda s} h(x_s) \,\diffd{s} + e^{-\lambda r} h(x_r)
            < \int_0^r \lambda e^{-\lambda s} h(x_s) \,\diffd{s} + e^{-\lambda r} V^{h,\pi}_\lambda(x_r).
    \end{equation}
    Hence, $\circled{\footnotesize1} < \circled{\footnotesize2}$ in \eqref{eq:app:recursive}, giving us that
    \begin{equation}
        \tilde{h}(x_0, \lambda) + e^{-\lambda 0} V^{h,\pi}_\lambda(x_0) = V^{h,\pi}_\lambda(x_0) = \sup_{0 \leq r \leq t} \left\{ \tilde{h}(x_r, \lambda) + e^{-\lambda r} V^{h,\pi}_\lambda(x_r) \right\}.
    \end{equation}
    Applying classical viscosity theory arguments and performing manipulations similar to \eqref{eq:app:supersol:p2} above yields
    \begin{equation}
        \dv{t} \left[
        \int_0^t \lambda e^{-\lambda \tau} h(x_\tau) \,\diffd{\tau} + e^{-\lambda t} V^{h,\pi}_{\lambda}(x_t)
        \right]_{t=0}
        = \nabla {V^{h,\pi}_{\lambda}}\T f(x_0, u_0) - \lambda( V^{h,\pi}_\lambda(x_0) - h(x_0))
        = 0 \geq 0.
    \end{equation}
    Hence, \eqref{eq:app:subsol:p1} also holds in this case.

    Combining both cases gives us that $V^{h,\pi}_\lambda$ is a viscosity sub-solution. Since $V^{h,\pi}_{\lambda}$ is both a viscosity super-solution and sub-solution to \eqref{eq:app:hj_ineq}, it is thus a viscosity solution of \eqref{eq:app:hj_ineq}.
\end{proof} \end{appendices}

\end{document}